\documentclass[11pt]{article}
\usepackage{latexsym}
\usepackage{amsfonts}
\usepackage{amssymb}
\usepackage{mathrsfs}

\newcommand{\be}{\begin{equation}}
\newcommand{\ee}{\end{equation}}
\newcommand{\bea}{\begin{eqnarray}}
\newcommand{\eea}{\end{eqnarray}}
\newcommand{\bean}{\begin{eqnarray*}}
\newcommand{\eean}{\end{eqnarray*}}
\newcommand{\brray}{\begin{array}}
\newcommand{\erray}{\end{array}}

\newcommand{\newsection}[1]{\setcounter{equation}{0}
\setcounter{dfn}{0}
\section{#1}}

\newtheorem{dfn}{Definition} [section]
\newtheorem{thm}[dfn]{Theorem}
\newtheorem{lmma}[dfn]{Lemma}
\newtheorem{ppsn}[dfn]{Proposition}
\newtheorem{crlre}[dfn]{Corollary}
\newtheorem{xmpl}[dfn]{Example}
\newtheorem{rmrk}[dfn]{Remark}

\newcommand{\bdfn}{\begin{dfn}}
\newcommand{\bthm}{\begin{thm}}
\newcommand{\blmma}{\begin{lmma}}
\newcommand{\bppsn}{\begin{ppsn}}
\newcommand{\bcrlre}{\begin{crlre}}
\newcommand{\bxmpl}{\begin{xmpl}}
\newcommand{\brmrk}{\begin{rmrk}}

\newcommand{\edfn}{\end{dfn}}
\newcommand{\ethm}{\end{thm}}
\newcommand{\elmma}{\end{lmma}}
\newcommand{\eppsn}{\end{ppsn}}
\newcommand{\ecrlre}{\end{crlre}}
\newcommand{\exmpl}{\end{xmpl}}
\newcommand{\ermrk}{\end{rmrk}}

\newcommand{\bbz}{\mathbb{Z}}

\newcommand{\cla}{\mathcal{A}}

\newcommand{\clh}{\mathcal{H}}

\newcommand{\clk}{\mathcal{K}}
\newcommand{\cll}{\mathcal{L}}

\newcommand{\clu}{\mathcal{U}}

\newcommand{\scrh}{\mathscr{H}}

\newcommand{\prf}{\noindent{\it Proof\/}: }

\newcommand{\nn}{\nonumber}

\def \qed { \mbox{}\hfill
$\Box$\vspace{1ex}}
\newcommand{\half}{\frac{1}{2}}
\newcommand{\id}{\mbox{id}}

\newcommand{\halpha}{\widehat{\alpha}}
\newcommand{\hbeta}{\widehat{\beta}}

\begin{document}

\author{{\sc Partha Sarathi Chakraborty}\hspace{-.2em}
and
{\sc Arupkumar Pal}}
\title{On equivariant Dirac operators for $SU_q(2)$
}
\maketitle

\begin{center}
\textit{Dedicated to Prof.\ Kalyan Sinha on his
sixtieth birthday.}
\end{center}

\begin{abstract}
We explain the notion of minimality for an equivariant spectral triple
and show that the triple for the quantum $SU(2)$ group 
constructed by Chakraborty and Pal in~(\cite{c-p1})
is minimal.
We also
give a decomposition of the spectral triple
constructed by Dabrowski et al~(\cite{dlssv})
in terms of the minimal  triple
constructed in~(\cite{c-p1}).
\end{abstract}
{\bf AMS Subject Classification No.:} {\large 58}B{\large 34}, {\large 46}L{\large 87}, {\large
  19}K{\large 33}

\section{Introduction}
The interaction between noncommutative geometry
and quantum groups, in particular the (noncommutative)
geometry of quantum groups, had been one of the less
understood and less explored areas of both the theories
for a while.
In the last few years, however, there has been some progress
in this direction.
The first important step was taken
by the authors in \cite{c-p1} where they found an optimal family
of Dirac operators for the quantum $SU(2)$ group acting on
$L_2(h)$, the $L_2$ space of the Haar state $h$, and equivariant
with respect to the (co-)action of the group itself.
This family has quite a few remarkable features, e.\ g.\
\begin{enumerate}
\item
any element
of the $K$-homology group can be realized by a member
from this family, which means that all elements of
the $K$-homology group are realizable through some Dirac operator
acting on the single Hilbert space $L_2(h)$ in a natural manner,
\item
the sign of any equivariant Dirac operator
on $L_2(h)$ is a compact perturbation of the sign of
a Dirac operator from this family,
\item
given any equivariant Dirac operator $\widetilde{D}$
acting on $L_2(h)$,
and any Dirac operator $D$ from this family,
there exist
two positive reals $k_1$ and $k_2$ such that
\[
|\widetilde{D}| \leq k_1 + k_2|D|,
\]
\item
they exhibit features that are unique
to the quantum case ($q\neq 1$).
It was proved in~\cite{c-p1} that for clasical $SU(2)$,
there does not exist any Dirac operator acting on (one copy of)
the $L_2$ space that is both equivariant as well as 3-summable.
\end{enumerate}
These triples were later analysed by Connes (\cite{co4})
in great detail, where the general theory of Connes-Moscovici
was applied to obtain a beautiful local index formula for
$SU_q(2)$.

Recently, Dabrowski et al (\cite{dlssv}) have constructed
another family of Dirac operators that act on
two copies of the $L_2$ space, has the right summability
property, is equivariant in a sense described in \cite{dlssv},
and is isospectral to the Classical Dirac operator.
In this note, we will give a decomposition of this Dirac
operator in terms of the Dirac operators constructed in
\cite{c-p1}.

\newsection{Equivariance and minimality}
In this section, we will formulate the notion of
an equivariant spectral triple for a compact quantum
group and what one means by its minimality, or irreducibility.
For basic notions on compact quantum groups, we refer the 
reader to~\cite{w3}. To fix notation, let us recall a few
things briefly here.
Let $G=(C(G),\Delta)$ be a compact quantum group, where
$C(G)$ is the unital $C^*$-algebra of `continuous functions on $G$'
and $\Delta$ the comultiplication map. The symbols $\kappa$ and $h$
will denote the antipode map and the Haar state for $G$.
For two functionals $\rho$ and $\sigma$ on $C(G)$, the convolution
product $\rho\ast\sigma$ is the functional $a\mapsto (\rho\otimes\sigma)\Delta(a)$.
For $\rho$ as above and $a\in C(G)$, we will denote by $a\ast\rho$ the
element $(\id\otimes\rho)\Delta(a)$ and by $\rho\ast a$ the element
$(\rho\otimes \id)\Delta(a)$.
A unitary representation $u$ of $G$ acting on a Hilbert space $\clh$
is   a unitary element of the multiplier algebra $M(\clk(\clh)\otimes C(G))$,
where $\clk(\clh)$ denotes the space of compact operators on $\clh$,
that satisfies the condition $(\id\otimes\Delta)u=u_{12}u_{13}$.
For a unitary representation $u$ and a continuous linear 
functional $\rho$ on $C(G)$, we will denote by $u_\rho$ the operator
$(\id\otimes\rho)u$ on $\clh$.
The GNS space associated with the state $h$ will be denoted by
$L_2(h)$ and the cyclic vector will be denoted by $\Omega$.
While using the comultiplication $\Delta$, we will often use
the Sweedler notation (i.e.\ $\Delta(a)=a_{(1)}\otimes a_{(2)}$).

Let $\cla$ be a unital $C^*$-algebra, $G$ be a
compact quantum group and let $\tau$ be an action
of $G$ on $\cla$, i.\ e.\ $\tau$ is a unital $C^*$-homomorphism
from $\cla$ to $\cla\otimes C(G)$ satisfying the condition
$(\id\otimes\Delta)\tau = (\tau\otimes\id)\tau$.
In other words, let $(\cla, G, \tau)$ be a $C^*$-dynamical system.
Recall (\cite{b-s}) that a covariant representation of $(\cla, G,\tau)$ on a Hilbert space
$\clh$ is a pair $(\pi,u)$ where $\pi$ is a unital *-representation
of $\cla$ on $\clh$, $u$ is a unitary representation of
$G$ on $\clh$
and they obey the condition
\be\label{cov_1}
u(\pi(a)\otimes I)u^* =(\pi\otimes\id)\tau(a),\quad a\in\cla.
\ee
By an \textbf{odd $G$-equivariant spectral data for $\cla$}, we mean
a quadruple $(\pi,u,\clh,D)$ where
\begin{enumerate}
\item
$(\pi,u)$ is a covariant representation of $(\cla, G,\tau)$ on 
the Hilbert space $\clh$,
\item
$\pi$ is faithful,
\item
$u(D\otimes I)u^* =D$,
\item
$(\pi,\clh,D)$ is an odd spectral triple.
\end{enumerate}
We will often be sloppy and just say $(\pi,\clh,D)$ is an odd $G$-equivariant spectral
triple for $\cla$, omitting $u$.
We say that an operator $D$ on a Hilbert space $\clh$ is an
\textbf{odd $G$-equivariant Dirac operator} for $\cla$ if there exists
a unitary representation $u$ of $G$ on $\clh$ such that
$(\pi,u,\clh,D)$ gives an $G$-equivariant spectral data for $\cla$.

Similarly, an even $G$-equivariant spectral data for $\cla$
consists of an even spectral data $(\pi,u,\clh, D,\gamma)$
where $(\pi,u,\clh,D)$ obeys conditions~1, 2 and~3 above, and moreover
$(\pi,\clh,D,\gamma)$ is an even spectral `triple', and one has
$u(\gamma\otimes I)u^* =\gamma$.
An even $G$-equivariant Dirac operator is also defined
similarly.

We say that an equivariant odd spectral data $(\pi,u,\clh,D)$ is \textbf{minimal}
if the covariant representation $(\pi,u)$ is irreducible.

Note that if we take $\cla=C(G)$, then the groups
$G$ and $G_{op}$ have natural actions $\Delta$ and $\Delta_{op}$
on $\cla$.
In what follows, we will mainly be concerned about
these two systems $(\cla=C(G), G, \Delta)$ and
$(\cla=C(G), G_{op}, \Delta_{op})$.
A $G$-equivariant spectral triple for $C(G)$ will be called
a \textbf{right equivariant spectral triple} for $C(G)$.
A \textbf{right equivariant Dirac operator} for $C(G)$ will mean a
$G$-equivariant Dirac operator for $C(G)$.
Similarly,
A $G_{op}$-equivariant spectral triple for $C(G)$ will be
called a \textbf{left equivariant spectral triple} for $C(G)$ and
a $G_{op}$-equivariant Dirac operator for $C(G)$
will be called a \textbf{left equivariant Dirac operator} for $C(G)$.

We will next study covariant representations of the
right $G$-action on $C(G)$, i.\ e.\ representations of the
system $(C(G), G,\Delta)$.
\blmma\label{tech_0}
Let $(\pi, u)$ be a covariant representation of
$(C(G), G, \Delta)$. If the Haar state $h$ of $G$ is faithful, 
then $\pi$ is faithful.
\elmma
\prf
Assume $\pi(a)=0$. Then $\pi(a^*a)=0$ and hence
$(\pi\otimes\id)\Delta(a^*a)=u(\pi(a^*a)\otimes I)u^*=0$.
Applying $(\id\otimes h)$ on both sides, we get
$h(a^*a)I=0$. Since $h$ is faithful, $a=0$.
\qed
\brmrk
{\rm
The above lemma helps ensure that if we have a compact quantum group
with a faithful Haar state, take a covariant representation $(\pi,u)$
of the system $(C(G), G,\Delta)$ on a Hilbert space $\clh$, and 
look at a Dirac operator $D$ on $\clh$, then we really get a 
spectral triple for the space $G$ rather than that of some
subspace (i.e.\ quotient $C^*$-algebra of $C(G)$) of it.
}
\ermrk

\blmma\label{tech_1}
Let $(\pi, u)$ be a covariant representation of
$(C(G), G, \Delta)$. Then the operator $u_h$ is a projection and for
any continuous linear functional $\rho$ on $\cla$,
one has $u_h u_\rho=u_\rho u_h=\rho(1)u_h$.
\elmma
\prf
Using Peter-Weyl decomposition for $u$, one can
assume without loss in generality that $u$ is finite dimensional.
Take two vectors $w$ and $w'$ in $\clh$.
Then
\bean
\langle w, u_h w'\rangle
 &=& \langle w\otimes\Omega , u ( w'\otimes\Omega)\rangle\\
 &=& \langle u^*(w\otimes\Omega) ,  w'\otimes\Omega\rangle\\
 &=& \langle ((\id\otimes\kappa)u)(w\otimes\Omega) ,
                 w'\otimes\Omega\rangle\\
 &=& \overline{\langle w'\otimes\Omega,
      ((\id\otimes\kappa)u)(w\otimes\Omega)\rangle}\\
 &=& \overline{\langle w', ((\id\otimes h\kappa)u) w\rangle}\\
 &=& \overline{\langle w', ((\id\otimes h)u) w\rangle}\\
 &=& \langle u_h w,  w'\rangle.
\eean
Thus $u_h$ is self adjoint.

Next, for any continuous linear functional $\rho$,
\bean
u_\rho u_h &=& (\id\otimes\rho)u (\id\otimes h)u\\
  &=& (\id\otimes\rho\otimes h)(u_{12}u_{13})\\
  &=& (\id\otimes\rho\otimes h)(\id\otimes\Delta)u\\
  &=& (\id\otimes \rho\ast h)u\\
  &=& \rho(1)u_h.
\eean
Similary one has $u_h u_\rho = \rho(1)u_h$.
In particular, $u_h^2=u_h$, so that $u_h$ is a projection.
\qed

\blmma\label{tech_2}
Let $A\equiv A(G)$ be the *-subalgebra of $C(G)$ generated
by matrix entries of all finite dimensional unitary
representations of $G$. Let $(A,\clu)$ be a dual pair of
Hopf$^*$-algebras (cf.\ \cite{v}). Then
\be\label{dlssv_eq1}
u_\rho \pi(a)=\pi(a\ast \rho_{(1)})u_{\rho_{(2)}}\quad
\mbox{for all }\rho\in\clu \mbox{ and } a\in A(G).
\ee
\elmma
\prf
Apply $(\id\otimes\rho)$ on both sides in the equality
\[
u(\pi(a)\otimes I)=\Bigl((\pi\otimes\id)\Delta(a)\Bigr)u
\]
and use the fact that $\rho(ab)=\rho_{(1)}(a)\rho_{(2)}(b)$.\qed

\blmma\label{tech_3}
Let $w\in\clh$ be a vector in the range of $u_h$.
Then for any $a\in A(G)$ and $\rho\in\clu$,
one has $u_\rho\pi(a)w=\pi(a\ast\rho)w$.
In particular, one has
$u_h \pi(a) w= h(a)w$.
\elmma
\prf
Use lemma~\ref{tech_2} to get
\bean
u_\rho\pi(a)w  &=& \pi(a\ast \rho_{(1)})u_{\rho_{(2)}}w\\
  &=& \pi(a\ast \rho_{(1)})u_{\rho_{(2)}} u_h w\\
  &=& \rho_{(2)}(1)\pi(a\ast \rho_{(1)}) u_h w\\
  &=& \pi(\rho_{(2)}(1)a\ast \rho_{(1)})  w\\
  &=& \pi(a_{(1)}\rho_{(1)}(a_{(2)})\rho_{(2)}(1))w \\
  &=& \pi(a_{(1)}\rho(a_{(2)}))w\\
  &=& \pi(a\ast \rho)w,
\eean
for $a\in A(G)$.
\qed

\blmma\label{tech_4}
The linear span of $\{\pi(a)u_h w: a\in A(G), w\in\clh\}$
is dense in $\clh$. In particular, $u_h$ is nonzero.
\elmma
\prf
Using Peter-Weyl  decomposition of $u$ and the observation
that $h(\kappa(a))=h(a)$ for all $a\in A$, it follows
that $u_h=(u^*)_h$.
Now take a vector $w'$ in $\clh$ such that
$\langle w', \pi(a)u_h w\rangle=0$ for all $w\in\clh$
and $a\in A$.
But then
$\langle w', \pi(a)(u^*)_h w\rangle=0$,
i.\ e.\
$\langle w'\otimes \Omega, (\pi(a)\otimes I)u^* (w\otimes \Omega)\rangle=0$.
The covariance condition~(\ref{cov_1}) now gives
$\langle u (w'\otimes\Omega), (\pi\otimes\id)\Delta(a) (w\otimes\Omega)\rangle=0$
for all $w\in\clh$ and $a\in A$.
In particular, one has
$\langle u (w'\otimes\Omega), (\pi\otimes\id)\Delta(a)
 (\pi(b)w\otimes\Omega)\rangle=0$
for all $w\in\clh$, and $a,b\in A$.
Since
$(\pi\otimes\id)\Delta(a)(\pi(b)\otimes I)=(\pi\otimes\id)(\Delta(a)(b\otimes I))$
and
$\{\Delta(a)(b\otimes I): a,b\in A\}$ is total
in $\cla\otimes\cla$, we get $u(w'\otimes\Omega)=0$
and consequently $w'=0$.\qed

For $w\in\clh$, denote by $P_w$ the projection onto
the closed linear span of $\{\pi(a)w: a\in A\}$.

\blmma\label{tech_5}
Let $w\in u_h\clh$.
Then
$(P_w\otimes I)u(P_w\otimes I)=u(P_w\otimes I)$.
If $w'$ is another vector in $u_h\clh$ such that $\langle w,w'\rangle=0$,
then the projections $P_w$ and $P_{w'}$ are orthogonal.
\elmma
\prf
For the first part, it is enough to show that
$P_w u_\rho P_w =u_\rho P_w$ for all $\rho\in\clu$.
But this is clear because from lemma~\ref{tech_3},
we have $u_\rho \pi(a) w = \pi(a\ast \rho) w$.

For the second part, take $a,a'\in A$.
Then using lemma~\ref{tech_4} one gets
\bean
\langle \pi(a)w,\pi(a')w'\rangle
  &=& \langle w, \pi(a^*a')w'\rangle\\
  &=& \langle u_h w, \pi(a^*a')w'\rangle\\
  &=& \langle w, u_h \pi(a^*a')w'\rangle\\
  &=& \langle w, h(a^*a')w'\rangle\\
  &=& 0.
\eean
Thus $P_w$ and $P_{w'}$ are orthogonal.
\qed

\bppsn\label{tech_6}
Let $\{w_1,w_2,\ldots\}$ be an orthonormal basis for
$u_h\clh$. Write $P_n$ for $P_{w_n}$, and let
$\pi_n(\cdot):=P_n \pi(\cdot)P_n$,
$u_n := (P_n\otimes I)u(P_n\otimes I)$.
Then
\begin{enumerate}
\item
for each $n$, $(\pi_n, u_n)$ is a covariant representation
of the system $(\cla, G,\Delta)$ on $P_n\clh$,
\item
$\pi=\oplus \pi_n$, $u=\oplus u_n$,
\item
$(\pi_n,u_n)$ is unitarily equivalent to the pair
$(\pi_L, u_R)$ where $\pi_L$ is the representation of $\cla$ on
$L_2(G)$ by left multiplications and $u_R$ is the right regular
representation of $G$.
\end{enumerate}
\eppsn
\prf
It follows from lemmas~\ref{tech_5} and \ref{tech_4} that
$P_n$'s are orthogonal, $\sum P_n=I$ and
consequently $\pi=\oplus \pi_n$ and $u=\oplus u_n$.

Define $V_n:P_n\clh\rightarrow L_2(G)$ by
\[
V_n \pi(a)w_n = \pi_L(a)\Omega,\quad a\in A.
\]
Since
$\langle \pi_L(a)\Omega,\pi_L(b)\Omega\rangle
 =h(a^*b)=\langle \pi(a)w_n,\pi(b) w_n\rangle$,
 $\{\pi(a)w_n:a\in A\}$ is total in $P_n\clh$
 and $\{\pi_L(a)\Omega:a\in A\}$ is total in $L_2(G)$,
 $V_n$ extends to a unitary from $P_n\clh$ onto $L_2(G)$.
 Next, for $a,b\in A$,
 one has
$V_n\pi(a)\pi(b)w_n=V_n\pi(ab)w_n=\pi_L(ab)\Omega
   =\pi_L(a)\pi_L(b)\Omega =\pi_L(a)V_n \pi(b)w_n$.
So $V_n \pi(a)=\pi_L(a)V_n$ for all $a\in A$ and hence for all $a\in\cla$.

Finally, we will show that $(V_n\otimes I)u(V_n^*\otimes I) = u_R$.
Write $\widetilde{u}_n:=(V_n\otimes I)u(V_n^*\otimes I)$.
Then  for any $\rho\in \clu$, one has
\bean
(\id\otimes \rho)\widetilde{u}_n \pi_L(a)\Omega
  &=& V_n u_\rho V_n^* \pi_L(a)\Omega\\
  &=& V_n u_\rho \pi(a) w_n \\
  &=& V_n u_\rho \pi(a) u_h w_n \\
  &=& V_n \pi(a\ast \rho) w_n\\
  &=& V_n \pi(a\ast \rho) V_n^* V_n w_n\\
  &=& \pi_L(a\ast\rho) \Omega.
\eean
By (\cite{w3}), $\widetilde{u}_n$ must be the right regular representation
$u$ on $L_2(G)$.
\qed

\brmrk
{\rm
The above proposition leads to an alternative proof
of the Takesaki-Takai duality for compact quantum groups
(theorem~7.5, \cite{b-s}).
}
\ermrk
\bthm
The covariant
representation $(\pi,u)$ is irreducible
if and only if the operator $u_h$ is a rank one projection.
\ethm
\prf
Immediate corollary of proposition~\ref{tech_6}.
\qed

\brmrk
In particular, it follows from the above theorem that
the covariant representation $(\pi_L, u_R)$
on $L_2(G)$ is irreducible. Thus the equivariant Dirac operator
constructed in~\cite{c-p1} is \textbf{minimal}.
\ermrk


\section{The decomposition}
\paragraph{Canonical triples for $SU_q(2)$.}
Let $q$ be a real number in the interval $(0,1)$.
Let $\cla$ denote the $C^*$-algebra of
continuous functions on $SU_q(2)$,
which is the universal $C^*$-algebra
generated by two elements $\alpha$ and $\beta$
subject to the relations
\[
\alpha^*\alpha+\beta^*\beta=I=\alpha\alpha^*+q^2\beta\beta^*,
\quad
\alpha\beta-q\beta\alpha=0=\alpha\beta^*-q\beta^*\alpha,
\quad
\beta^*\beta=\beta\beta^*
\]
as in \cite{c-p1}.
Let $\pi :\cla\rightarrow\cll(L_2(h))$ be
the representation given by left multiplication
by elements in $\cla$.
Let $u$ denote the right regular representation of
$SU_q(2)$.
Recall~(\cite{w3}) that $u$ is the unique representation acting on
$L_2(h)$ that obeys the condition
\be
\Bigl((\mbox{id}\otimes\rho)u\Bigr)\pi(a)\Omega
 = \pi\Bigl((\mbox{id}\otimes\rho)\Delta(a)\Bigr)\Omega
\ee
for all $a\in\cla$ and
for all continuous linear functionals $\rho$ on $\cla$.
In \cite{c-p1}, the authors studied right equivariant Dirac operators,
those Dirac operators that commute with the right regular representation,
i.\ e.\  $D$ acting on $L_2(h)$ for which
\[
(D\otimes I)u=u(D\otimes I).
\]
In particular, an optimal family of equivariant Dirac
operators were found. A generic member of this family
is of the form
\[
e^{(n)}_{ij}\mapsto\cases{ (an+b)e^{(n)}_{ij} & if $-n\leq i< n-k$,\cr
        (cn+d)e^{(n)}_{ij} & if $i=n-k,n-k+1,\ldots,n$,}
\]
where $k$ is a fixed nonnegative integer
and $a$, $b$, $c$, $d$ are reals with $ac<0$.
If one looks at left equivariant Dirac operators,
the same arguments would then lead to the following
theorem.
\bthm
Let $v$ be the left regular representation of $SU_q(2)$.
Let $k$ be a nonnegative integer and let $a$, $b$, $c$, $d$ be real
numbers with $ac<0$.
Then the operator $D\equiv D(k,a,b,c,d)$ on $L_2(h)$ given by
\[
e^{(n)}_{ij}\mapsto
   \cases{ (an+b)e^{(n)}_{ij} & if $-n\leq j< n-k$,\cr
         (cn+d)e^{(n)}_{ij} & if $j=n-k,n-k+1,\ldots,n$,}
\]
gives a spectral triple $(\pi,L_2(h),D)$ having nontrivial
Chern character and obeys
\be\label{equiv}
(D\otimes I)v=v(D\otimes I).
\ee

Conversely, given any spectral triple $(\pi,L_2(h), \widetilde{D})$
with nontrivial Chern character such that
$(\widetilde{D}\otimes I)v=v(\widetilde{D}\otimes I)$,
there exist a nonnegative integer $k$ and reals $a$, $b$, $c$, $d$ with $ac<0$
such that
\begin{enumerate}
\item
$\mbox{sign\,}\widetilde{D}$ is a compact perturbation of
the sign of $D\equiv D(k,a,b,c,d)$, and
\item
there exist constants $k_1$ and $k_2$ such that
\[
 |\widetilde{D}|\leq k_1+k_2|D|.
 \]
 \end{enumerate}
\ethm
\prf
The key point is to note that
the characterizing property
of the left regular representation $v$ is
\be
\Bigl((\mbox{id}\otimes\rho)v^*\Bigr)\pi(a)\Omega
 = \pi\Bigl((\rho\otimes\mbox{id})\Delta(a)\Bigr)\Omega.
\ee
Thus on the right hand side, one now has left convolution
of $a$ by $\rho$ instead of right convolution by $\rho$.
Therefore any self-adjoint operator on $L_2(h)$ with discrete
spectrum that obeys $(D\otimes I)v=v(D\otimes I)$ will be of the form
\[
e^{(n)}_{ij}\mapsto \lambda(n,j)e^{(n)}_{ij}.
\]
Hence if one now proceeds exactly along the
same lines as in \cite{c-p1},
one gets all the desired conclusions.
\qed

Observe at this point that the whole analysis carried out in~\cite{co4}
will go through for this Dirac operator as well.
Let us now take two such Dirac operators $D_1$ and $D_2$ on $L_2(h)$
given by
\be
D_1 e^{(n)}_{ij}=\cases{-2n e^{(n)}_{ij} & if $j\neq n$,\cr
                       (2n+1) e^{(n)}_{ij} & if $j=n$,}
\qquad
D_2 e^{(n)}_{ij}=\cases{(-2n-1) e^{(n)}_{ij} & if $j\neq n$,\cr
                       (2n+1) e^{(n)}_{ij} & if $j=n$.}
\ee
Now look at the triple
\[
(\pi\oplus\pi, L_2(h)\oplus L_2(h), D_1\oplus |D_2|).
\]
It is easy to see that this is a spectral triple.
Nontriviality of its Chern character is a direct
consequence of that of $D_1$.
We will show in the next paragraph that in a certain sense,
the spectral triple constructed in \cite{dlssv} is equivalent
to this above triple.

\paragraph{The decomposition.}
Let us briefly recall the Dirac operator
constructed in~\cite{dlssv}.
The carrier Hilbert space $\scrh$ is a direct sum of two
copies of $L_2(h)$ that decomposes as
\[
\scrh=W_0^\uparrow\oplus \Bigl(\bigoplus_{n\in\half\bbz_+}
    ( W_n^\uparrow\oplus W_n^\downarrow )\Bigr),
\]
where
\bean
W_n^\uparrow &=&\mbox{span}\{ u^n_{ij}:i=-n,-n+1,\ldots,n,\,
     j=-n-\half,-n+\half,\ldots,n+\half\},\\
W_n^\downarrow &=&\mbox{span}\{ d^n_{ij}:i=-n,-n+1,\ldots,n,\,
     j=-n+\half,-n+\frac{3}{2},\ldots,n-\half\}.
\eean
($u^n_{ij}$ and $d^n_{ij}$ correspond to the basis elements
$|nij\!\uparrow\rangle$ and $|nij\!\downarrow\rangle$
respectively in the notation of~\cite{dlssv}).
Now write
\[
v^n_{ij}=\left(\matrix{u^n_{ij} \cr d^n_{ij}}\right)
\]
with the convention that $d^n_{ij}=0$ for $j=\pm(n+\half)$.
Then the  representation $\pi'$ of $\cla$ on
$\scrh$ is given by
\bean
\pi'(\alpha^*) \,v^n_{ij}
&=&  a^+_{nij} v^{n+\half}_{ i+\half, j+\half}
 +  a^-_{nij} v^{n-\half}_{i+\half, j+\half},
\nn \\
\pi'(-\beta) \,v^n_{ij}
&= & b^+_{nij} v^{n+\half}_{ i+\half, j-\half}
 +  b^-_{nij} v^{n-\half}_{ i+\half, j-\half},
\nn \\
\pi'(\alpha) \,v^n_{ij}
&= &\tilde a^+_{nij} v^{n+\half}_{ i-\half, j-\half}
 + \tilde a^-_{nij} v^{n-\half}_{ i-\half, j-\half},
 \\
\pi'(-\beta^*) \, v^n_{ij}
&= &\tilde b^+_{nij} v^{n+\half}_{ i-\half, j+\half}
 + \tilde b^-_{nij} v^{n-\half}_{ i-\half, j+\half},
\nn
\eean
where $ a^\pm_{nij}$ and $ b^\pm_{nij}$ are
the following  $2 \times 2$ matrices:
\bean
 a^+_{nij} &=& q^{(i+j-\half)/2} [n + i + 1]^\half
\left(\matrix{
q^{-n-\half} \, \frac{[n+j+\frac{3}{2}]^{1/2} }{[2n+2]} & 0 \cr
q^\half \,\frac{[n-j+\half]^{1/2}}{[2n+1]\,[2n+2]} &
q^{-n} \, \frac{[n+j+\half]^{1/2}}{[2n+1]}
}\right),
\nn \\[3ex]
 a^-_{nij} &= &q^{(i+j-\half)/2} [n - i]^\half
\left(\matrix{
q^{n+1} \, \frac{[n-j+\half]^{1/2}}{[2n+1]} &
- q^\half \,\frac{[n+j+\half]^{1/2}}{[2n]\,[2n+1]} \cr
0 & q^{n+\half} \, \frac{[n-j-\half]^{1/2}}{[2n]}
}\right),
\nn \\[3ex]
 b^+_{nij} &=& q^{(i+j-\half)/2} [n + i + 1]^\half
\left(\matrix{
\frac{[n-j+\frac{3}{2}]^{1/2}}{[2n+2]} & 0 \cr
- q^{-n-1} \,\frac{[n+j+\half]^{1/2}}{[2n+1]\,[2n+2]} &
q^{-\half} \, \frac{[n-j+\half]^{1/2}}{[2n+1]}
}\right),
\\[3ex]
 b^-_{nij} &=& q^{(i+j-\half)/2} [n - i]^\half
\left(\matrix{
- q^{-\half} \, \frac{[n+j+\half]^{1/2}}{[2n+1]} &
- q^n \,\frac{[n-j+\half]^{1/2}}{[2n]\,[2n+1]} \cr
0 & - \frac{[n+j-\half]^{1/2}}{[2n]}
}\right),
\nn
\eean
($[m]$ being the $q$-number $\frac{q^m-q^{-m}}{q-q^{-1}}$)
and $\tilde a^\pm_{nij}$  and $\tilde b^\pm_{nij}$
are the hermitian conjugates of the above ones:
\[
\tilde a^\pm_{nij} = (a^\mp_{n\pm\half, i-\half, j-\half})^*,
\qquad
\tilde b^\pm_{nij} = (b^\mp_{n\pm\half, i-\half, j+\half})^*.
\]
The operator $D$ is given by
\[
D u^n_{ij}=(2n+1) u^n_{ij},
\qquad
D d^n_{ij}= -2n d^n_{ij}.
\]
The triple $(\pi', \scrh, D)$ is precisely the triple
constructed in~\cite{dlssv}.

\bthm
Let $\clk_q$ be the two-sided ideal of $\cll(\mathscr{H})$
generated by the operator
\[
d^n_{ij}\mapsto q^n d^n_{ij},
\quad
u^n_{ij}\mapsto q^n u^n_{ij},
\]
and let $\cla_f$ denote the *-subalgebra of $\cla$ generated
by $\alpha$ and $\beta$.
Then there is a unitary $U:L_2(h)\oplus L_2(h)\rightarrow \scrh$
such that
\bea
U(D_1\oplus |D_2|)U^*  &= &  D,\\
U\Bigl(\pi(a)\oplus\pi(a)\Bigr)U^* -\pi'(a)
     &\in & \clk_q \qquad \mbox{for all } a\in\cla_f.
\eea
\ethm
\prf
Define $U:L_2(h)\oplus L_2(h)\rightarrow \scrh$
as follows
\bean
U(e^{(n)}_{ij}\oplus 0) &=& d^n_{i,j+\half},
\quad
i=-n,-n+1,\ldots,n,\;
j=-n,-n+1,\ldots,n-1,\\
U(e^{(n)}_{in}\oplus 0) &=& u^n_{i,n+\half},
\quad
i=-n,-n+1,\ldots,n,\\
U(0\oplus e^{(n)}_{ij}) &=& u^n_{i,j-\half},
\quad
i=-n,-n+1,\ldots,n,\;
j=-n,-n+1,\ldots,n.
\eean
It is immediate that $U(D_1\oplus |D_2|)U^* = D$.
Therefore all that we need to prove now is that
$U(\pi(a)\oplus\pi(a))U^* -\pi'(a)\in\clk_q$
for all $a\in \cla_f$.
For this, let us introduce the representation
$\widehat{\pi}:\cla\rightarrow\cll(L_2(h))$
given by
\[
\widehat{\pi}(\alpha)=\halpha,
\qquad
\widehat{\pi}(\beta)=\hbeta,
\]
where $\halpha$ and $\hbeta$ are the following operators
on $L_2(h)$ (see lemma~2.2, \cite{c-p3})
\bea
\halpha: e^{(n)}_{ij}  &\mapsto&
q^{2n+i+j+1} e^{(n+\half)}_{i-\half ,j-\half }
    + (1-q^{2n+2i})^\half(1-q^{2n+2j})^\half
    e^{(n-\half )}_{i-\half ,j-\half },\label{halpha}\\
\hbeta:e^{(n)}_{ij}  &\mapsto&
   - q^{n+j}(1-q^{2n+2i+2})^\half
            e^{(n+\half )}_{i+\half ,j-\half }
  + q^{n+i}(1-q^{2n+2j})^\half
 e^{(n-\half )}_{i+\half ,j-\half },\label{hbeta}
\eea
It is easy to see that
\[
\pi(a)\oplus\pi(a)-\widehat{\pi}(a)\oplus\widehat{\pi}(a)
\in U^*\clk_q U.
\]
for $a=\alpha^*$ and $a=\beta$.
Therefore it is enough to verify that
\[
U(\widehat{\pi}(a)\oplus\widehat{\pi}(a))U^*-\pi'(a)\in\clk_q
\]
for $a=\alpha^*$ and for $a=\beta$.

Next observe that
\bean
 a^+_{nij} &=& (1-q^{2n+2i+2})^\half
\left(\matrix{
 (1-q^{2n+2j+3})^\half  & 0 \cr
0 &
 (1-q^{2n+2j+1})^\half
}\right) + O(q^{2n}),
 \\[3ex]
 a^-_{nij} &= &q^{2n+i+j+\half}(1-q^{2n-2i})^\half
\left(\matrix{
q(1-q^{2n-2j+1})^\half  & 0 \cr
0 & (1-q^{2n-2j-1})^\half
}\right)+O(q^{2n}),
 \\[3ex]
 b^+_{nij} &=& q^{n+j-\half}(1-q^{2n+2i+2})^\half
\left(\matrix{
q(1-q^{2n-2j+3})^\half  & 0 \cr
0 & (1-q^{2n-2j+1})^\half
}\right)+O(q^{2n}),
\\[3ex]
&=& q^{n+j-\half}(1-q^{2n+2i+2})^\half
\left(\matrix{
q  & 0 \cr
0 & 1
}\right)+O(q^{2n}),
\\[3ex]
 b^-_{nij} &=&- q^{n+i}(1-q^{2n-2i})^\half
\left(\matrix{
(1-q^{2n+2j+1})^\half  & 0
 \cr
0 &(1-q^{2n+2j-1})^\half
}\right)+O(q^{2n})\\[3ex]
&=& - q^{n+i}
\left(\matrix{
(1-q^{2n+2j+1})^\half  & 0
 \cr
0 &(1-q^{2n+2j-1})^\half
}\right)+O(q^{2n}).
\eean
The required result now follows from this easily.
\qed

\brmrk
\emph{
The above decomposition in particular tells us that the spectral triples
$(\pi\oplus\pi, L_2(h)\oplus L_2(h), D_1\oplus |D_2|)$
and
$(\pi', \scrh, D)$
are essentially unitarily equivalent at the Fredholm module level.
therefore by proposition~8.3.14, \cite{h-r}, they give rise to the same element
in K-homology.
}
\ermrk

\brmrk
{\rm
In the spectral triple in~\cite{dlssv},  the Hilbert space can be
decomposed as a direct sum of two isomorphic copies in such a manner 
that in each half Dirac operator
has constant sign. So positive and negative signs come with equal
frequency. However this symmetry is only superficial, as the 
decomposition above illustrates.
This asymmetry might be a reflection of the inherent asymmetry in the
\textit{growth graph} associated with quantum $SU(2)$ (cf.\ \cite{c-p}).
For classical $SU(2)$ the graph is symmetric whereas in
the quantum case it is not.

It should also be pointed out here that, at least as far as classical
odd dimensional spaces are concerned, this kind of sign symmetry is
always superficial. They are always inherent  in the
even cases, not in the odd ones.
}
\ermrk

{\footnotesize {\bf Acknowledgements.} 
 We would like to thank the referee for suggesting some improvements.}


\noindent{\sc Partha Sarathi Chakraborty}
(\texttt{parthac@imsc.res.in})\\
         {\footnotesize  Institute of Mathematical Sciences, 
   CIT Campus, Taramani, Chennai--600\,113, INDIA}\\[1ex]
{\sc Arupkumar Pal} (\texttt{arup@isid.ac.in})\\
         {\footnotesize Indian Statistical
Institute, 7, SJSS Marg, New Delhi--110\,016, INDIA}

\end{document}